\documentclass[11pt]{amsart}
\usepackage{amsmath,amssymb,amsthm}
\usepackage[a4paper,margin=28mm]{geometry}
\usepackage[hidelinks]{hyperref}
\newtheorem{mainthm}{Theorem}

\newtheorem{theorem}{Theorem}[section]
\newtheorem{proposition}[theorem]{Proposition}
\newtheorem{lemma}[theorem]{Lemma}
\newtheorem{corollary}[theorem]{Corollary}
\theoremstyle{remark}
\newtheorem{example}[theorem]{Example}

\newcommand{\F}{\mathbb F_p}
\newcommand{\Inn}{\operatorname{Inn}}
\newcommand{\Aut}{\operatorname{Aut}}
\newcommand{\Hom}{\operatorname{Hom}}

\newcommand{\Der}{\operatorname{Der}}
\title[Berkovich's conjecture]{Berkovich's conjecture on noninner automorphisms of order $p$}
\author{M. Yas\.{I}r K{\i}zmaz}
\address{Department of Mathematics, Bilkent University, Ankara, Turkey}
\begin{document}
\begin{abstract}
We prove a long-standing conjecture of Berkovich (Kourovka Notebook,
Problem~4.13): every finite nonabelian $p$-group $P$ has a noninner
automorphism of order $p$. If $p$ is odd, the
automorphism $\alpha$ can be chosen with $[P,\alpha,\alpha]=1$. More
precisely, given a maximal normal elementary abelian subgroup $E$ of
$P$, we can choose $\alpha$ to fix both $F$ and $P/F$ pointwise, where
$F$ is either $E$ or a hyperplane of $E$ normal in $P$. If $p=2$, the
automorphism can be chosen with $[P,\alpha,\alpha,\alpha]=1$, and with
$[P,\alpha,\alpha]=1$ when $Z(P)$ is noncyclic. The dihedral group of
order eight shows that the hypothesis on $Z(P)$ cannot be omitted.
\end{abstract}
\maketitle

\section{Introduction}

Let $P$ be a finite nonabelian $p$-group. By a 1966 theorem of
Gaschütz, $P$ has a noninner automorphism $\alpha$ of $p$-power
order \cite{Gaschutz}. Assume that $\alpha$ has order $p^n$. Note that
$n\geq 1$ as $\alpha$ is not inner. Then clearly the automorphism
$\alpha^{p^{n-1}}$ has order $p$, but it may be inner.
In 1973, Berkovich asked whether every finite nonabelian $p$-group has a
noninner automorphism of order $p$ \cite[Problem~4.13]{Kourovka}.
This question is now known as Berkovich's conjecture. Although it is
simple to state, it has remained open for more than fifty years. Many
partial results towards it are collected in the survey \cite{JS}.
In this article, we prove the conjecture for all primes, with no
additional hypothesis on $P$.

For an automorphism $\alpha$ of $P$ and an $\alpha$-invariant subgroup
$H$ of $P$, we write $[H,\alpha]=\langle h^{-1}\alpha(h):h\in H\rangle$,
and we take iterated commutators from the left. We say that $\alpha$
is \emph{quadratic} if $[P,\alpha,\alpha]=1$, and \emph{cubic} if
$[P,\alpha,\alpha,\alpha]=1$. The proofs for odd $p$ and for $p=2$
share the general tools but differ in the final steps, so we state the
two cases separately. Our first theorem treats odd primes.

\begin{mainthm}\label{berkovich}
Let $p$ be an odd prime. Then every finite nonabelian $p$-group $P$
admits a noninner automorphism $\alpha$ of order $p$ with
$[P,\alpha,\alpha]=1$.
\end{mainthm}

For $p=2$, we obtain a similar result with a weaker commutator condition.

\begin{mainthm}\label{two}
Every finite nonabelian $2$-group $P$ admits a noninner automorphism
$\alpha$ of order two with $[P,\alpha,\alpha,\alpha]=1$. If $Z(P)$ is
noncyclic, then $\alpha$ can be chosen with $[P,\alpha,\alpha]=1$.
\end{mainthm}

The commutator condition in Theorem~\ref{two} is sharp, since every
noninner involution $\beta$ of the dihedral group $D_8$ of order eight
satisfies $[D_8,\beta,\beta]\ne1$ (Example~\ref{d8}). Theorems~\ref{berkovich}
and~\ref{two} together give the following.

\begin{corollary}\label{conjecture}
Let $p$ be a prime. Then every finite nonabelian $p$-group has a
noninner automorphism of order $p$. In particular, Berkovich's
conjecture holds.
\end{corollary}

For $2$-groups, Abdollahi and Ghoraishi obtained a noninner
automorphism of order two or four fixing $\Phi(P)$ pointwise when $P$
has coclass two or $(P,Z(P))$ is a Camina pair \cite{AGtwo}.

Berkovich's conjecture has been proven for several classes of groups.
Liebeck proved the odd-prime case for groups of class two
\cite{Liebeck}, and Abdollahi treated class two for all primes
\cite{Abdollahi}. He later extended this to groups $P$ for which
$P/Z(P)$ is powerful \cite{Abdpow}. Jamali and Viseh proved the conjecture for
groups with cyclic derived subgroup \cite{JV}. The result for groups
of class three is due to
Abdollahi, Ghoraishi and Wilkens \cite{AGW}. For $p$-groups of
coclass two, Abdollahi, Ghoraishi, Guerboussa, Reguiat and Wilkens
proved that there is a noninner automorphism of order $p$ fixing
$Z(P)$ pointwise \cite{AGGRW}. Ruscitti, Legarreta and Yadav proved the
conjecture for groups of coclass three when $p\ne3$ \cite{RLY}.
Komma proved the conjecture for groups of coclass four or five when $p\ge5$
\cite{Komma}. A cohomological approach goes back to Schmid, who
proved that if $P$ is regular, $1\ne N\trianglelefteq P$ and $P/N$
is noncyclic, then $\widehat H^n(P/N,Z(N))\ne0$ for every integer
$n$ \cite{Schmid}. Benmoussa and Guerboussa proved the same
nonvanishing result for semi-abelian $p$-groups under the additional
hypothesis that $P/N$ is not generalized quaternion. They also proved
Berkovich's conjecture for semi-abelian $p$-groups
\cite[Theorems~1.1 and~1.3]{BG}.

Deaconescu and Silberberg proved that if
$C_P(Z(\Phi(P)))\ne\Phi(P)$, then $P$ has a noninner automorphism
of order $p$ fixing $\Phi(P)$ pointwise \cite{DS}. This reduces the
conjecture to the case $C_P(Z(\Phi(P)))=\Phi(P)$. Let $W_0$ be the
preimage of $\Omega_1(Z_2(P)/Z(P))$ in $Z_2(P)$. Ghoraishi
strengthened this reduction by showing that such an automorphism
exists unless $W_0\le C_P(W_0)=\Phi(P)$
\cite[Theorem~1.1]{Ghoraishi}. Singh, Garg and Kalra later gave a
short elementary proof of Ghoraishi's result \cite[Theorem~2.6]{SGK}.

Garg and Singh treated the case where $Z(P)$ is noncyclic and
$W_0$ is nonabelian \cite{GS}. This case follows from Ghoraishi's
theorem because $W_0\le C_P(W_0)$ forces $W_0$ to be abelian.
Singh and Sharma proved the conjecture for groups of odd order
with cyclic center under additional hypotheses \cite{SS}.
For odd $p$ and cyclic center, Ma and Xu proved it when
$C_P(Z(\gamma_3(P)P^p))\nleq\gamma_3(P)P^p$ \cite{MaXu}.

In a preprint, Xu proposed a proof for all
primes \cite{Xu}. As far as we can see, some steps of that argument do
not hold as stated, and the problem is still listed as open in the
latest edition of the Kourovka Notebook \cite{Kourovka}.

We explain how we choose the automorphism in
Theorem~\ref{berkovich}. Given a normal elementary abelian subgroup
$F$ of $P$, put
$$
 D_F(P)=\{\alpha\in\Aut(P):\alpha|_F=1,\
                       \alpha\text{ induces the identity on }P/F\}.
$$
If $\alpha\in D_F(P)$ and $g\in P$, then $\alpha(g)=gf_g$ for some
$f_g\in F$. Since $\alpha$ fixes $f_g$, we have $\alpha^p(g)=gf_g^p=g$.
Thus every nonidentity element of $D_F(P)$ has order $p$. Moreover,
$[P,\alpha]\le F$ and $[F,\alpha]=1$, so every element of $D_F(P)$ is
quadratic. The question is how to choose $F$ so that some element of
$D_F(P)$ is noninner.

Let $E$ be maximal among the normal elementary abelian subgroups of $P$.
If every element of $D_E(P)$ is inner, we say that $E$ is
\emph{saturated}. We use the same terminology for any normal elementary
abelian subgroup $F$ satisfying $D_F(P)\le\Inn(P)$. The base group of
$C_p\wr C_p$ gives an example of a saturated $E$
(Example~\ref{wreath}). Our next theorem shows that for odd $p$ it
suffices in this situation to pass to a hyperplane of $E$.

\begin{mainthm}\label{flag}
Let $p$ be odd, let $P$ be a finite $p$-group which is not elementary
abelian, and let $E$ be maximal by inclusion among its normal elementary
abelian subgroups. Put $V=\Omega_1(Z(P))$ and assume that $E$ is
saturated. Then $V<E$, and for every $P$-invariant hyperplane $F$
of $E$ there is a noninner automorphism $\alpha\in D_F(P)$ of order $p$ such that
$\alpha|_E=1+s$, where $\ker s=F$ and $s(E)$ is a line in $V$.
\end{mainthm}

In particular, for odd $p$, every chain $V=F_0<F_1<\cdots<F_m=E$ of normal
subgroups of $P$ with $|F_i:F_{i-1}|=p$ contains an unsaturated term.
If $E$ is saturated, then $m\ge1$ and $F_{m-1}$ is unsaturated by
Theorem~\ref{flag}. Otherwise $F_m=E$ is unsaturated. Since every
element of $D_F(P)$ is quadratic, Theorem~\ref{flag} and the
unsaturated case give Theorem~\ref{berkovich}.

Sections~\ref{sec:cosets}--\ref{sec:lift} are written for an arbitrary
prime, except for Lemma~\ref{regular}\textup{(iii)},
Example~\ref{wreath}, Corollary~\ref{noncyclic-centre} and
Lemma~\ref{inner-lift}, which assume that $p$ is odd. For odd $p$, when $E$ is saturated, we lift a transvection
$1+s$ of $E$, with $s(E)\le V$ and $s([E,P])=0$, to an automorphism
inducing the identity on $P/E$. The lift has order $p$ by
Lemma~\ref{iterate}.
If $C_P(E)>E$, Proposition~\ref{structure} shows that $C_P(E)/E$ is
cyclic and that its unique subgroup of order $p$ is central in $P/E$.
We first realize $s$ as the restriction of a derivation on the
preimage of this subgroup in $P$. A dimension count lets us apply
Lemma~\ref{extend} to extend $s$ to an $E$-valued derivation on $P$
(Lemma~\ref{lift-transvection}). We then use transfer and a normal
subgroup of order $p^2$ in $P/E$ to show that every such lift induces
the identity on $P/[E,P]$ (Proposition~\ref{displacement}). In
particular, the lift belongs to $D_F(P)$ for $F=\ker s$.

For $p=2$, the elements of $D_E(P)$ are again quadratic of order at
most two, so the main case is again the one in which $E$ is saturated.
Put $Y=\{y\in C_P(E):[P,y]\le E\}$. When $E$ is saturated, every class
in $H^1(P/E,E)$ is induced by conjugation by an element of $Y$. If
$\dim H^1(P/E,E)>\dim\Omega_1(Z(P))$, then $Y$ contains a quaternion
subgroup of order eight. Inverting a cyclic subgroup of order four of this quaternion
group gives the required involution (Proposition~\ref{invertcase}).
Otherwise, either $P$ splits over $E$, or both $Z(P)$ and $C_P(E)/E$
are cyclic. In the latter case, restriction to the subgroup of order
two of $C_P(E)/E$ detects the extension class of $P$, and
Lemma~\ref{lift-transvection} lifts a transvection of $E$. If $P/E$
has a normal four-subgroup, the lift can be chosen to be a quadratic
involution. Otherwise $P/E$ is dihedral, semidihedral or generalized
quaternion. Here the two cases of the proof diverge. For a lift
$\alpha(g)=g\delta(g)$ of $1+s$, Lemma~\ref{iterate} gives $\alpha^p=1$
for odd $p$, but only $\alpha^2(g)=g\,s(\delta(g))$ for $p=2$. In the
last part of the proof, we compute this square, and we correct the lift
by a derivation with values in a cyclic subgroup of order four so that
the corrected automorphism has order two.

\section{Derivations on cosets}\label{sec:cosets}

All modules in this section are finite-dimensional right modules over
$k=\mathbb F_p$. We write $\Der(L,M)$ for the space of maps satisfying
$d(gh)=d(g)^h+d(h)$. The principal derivations have the form
$g\mapsto m^g-m$. Put
\[
 b_L(M)=\dim\Der(L,M)-\dim M,\qquad
 h_L(M)=b_L(M)+\dim M^L.
\]
Thus $h_L(M)$ is the dimension of the space of derivations modulo
principal derivations.

Let $H\le L$, put $X=L/H$, and write $T=\operatorname{Map}(X,M)$ for
the space of functions on the left cosets, with $f^g(x)=f(gx)^g$.
We use the inclusion $\eta:M\to T$ by constant functions and the
map $\varepsilon:T\to M$ given by summing the values.
For each $x\in X$, we choose $r_x\in x$ with $r_H=1$.
If $d:H\to M$ is a derivation, put
\begin{equation}\label{cosetformula}
 \Delta(g)(x)=d(r_{gx}^{-1}g r_x)^{r_x^{-1}}.
\end{equation}
Set $n(g,x)=r_{gx}^{-1}g r_x$. Then $n(gh,x)=n(g,hx)n(h,x)$, and
substitution shows that $\Delta(gh)=\Delta(g)^h+\Delta(h)$. Moreover,
$\Delta(u)(H)=d(u)$ for $u\in H$.

\begin{lemma}\label{restriction}
With the notation above, assume that $L$ is a finite $p$-group and
put $m=[L:H]$.
\begin{enumerate}
\item[\textup{(i)}] Evaluation $\Delta\mapsto(u\mapsto \Delta(u)(H))$ maps
$\Der(L,T)$ onto $\Der(H,M)$. Its kernel consists of the principal
derivations $g\mapsto f^g-f$ with $f(H)=0$.
\item[\textup{(ii)}] We have $b_H(M)\le m b_L(M)$. In case of
equality, both maps
\[
 \Der(L,T)\longrightarrow\Der(L,T/\eta M),\qquad
 \Der(L,T)\xrightarrow{\ \varepsilon\ }\Der(L,M)
\]
are onto.
\end{enumerate}
\end{lemma}
\begin{proof}
Formula~\eqref{cosetformula} gives surjectivity in (i). If $f(H)=0$,
then $(f^u-f)(H)=0$ for every $u\in H$. Conversely, let $\Delta$ be in
the kernel and put $f(x)=\Delta(r_x)(H)^{r_x^{-1}}$. Evaluating
$\Delta(gr_x)=\Delta(r_{gx}n(g,x))$ at $H$, we obtain that
$\Delta(g)=f^g-f$. An invariant function vanishing at $H$ vanishes
everywhere. Thus the kernel has dimension $(m-1)\dim M$, and so
$b_L(T)=b_H(M)$.

For any submodule $A$ of $B$ we have
$b_L(B)\le b_L(A)+b_L(B/A)$, with equality precisely when every
derivation into $B/A$ lifts to $B$. This follows by taking the kernel
of $\Der(L,B)\to\Der(L,B/A)$. The map $\varphi\otimes v\mapsto
(x\mapsto\varphi(x)v)$ identifies $\operatorname{Map}(X,k)\otimes_k M$,
with the diagonal action, with $T$. Every composition factor of
$\operatorname{Map}(X,k)$ is trivial since $L$ is a $p$-group. We
tensor a composition series of $\operatorname{Map}(X,k)$ with $M$. This yields a filtration of $T$ with $m$ factors
isomorphic to $M$, and so $b_H(M)=b_L(T)\le m\,b_L(M)$.

Now assume that equality holds. We may start the series with the
constant line, or end it with the quotient given by summation.
Applying the same bound to the other $m-1$ factors, we see that the
inequality for the first (respectively last) step is an equality.
Hence every derivation into $T/\eta M$ (respectively into $M$) lifts
to $T$, as asserted.
\end{proof}

In cohomological terms, $b_L(M)=\dim H^1(L,M)-\dim H^0(L,M)$.
Under Shapiro's isomorphism, the maps induced by $\eta$ and
$\varepsilon$ are restriction and corestriction, respectively
\cite[III.6 and III.9]{Brown}. The second surjectivity in (ii) is
equivalent to surjectivity of $H^1(H,M)\to H^1(L,M)$, since a
principal derivation into $M$ lifts to $T$ as $\varepsilon$ is onto.
For $M=k$ with trivial action, (ii) gives the Schreier inequality
$d(H)-1\le [L:H](d(L)-1)$, where $d(L)$ is the minimal number of
generators of $L$.

\begin{lemma}\label{extend}
Let $K\trianglelefteq R$, let $L=R/K$ be a finite $p$-group, and
let $M$ be a $kL$-module. Let $J$ be the inverse image of $H\le L$.
Assume that $b_H(M)=[L:H]b_L(M)$. If an $R$-equivariant map
$s:K\to M$ extends to a derivation on $J$, it extends to a
derivation on $R$.
\end{lemma}
\begin{proof}
We identify $R/J$ with $L/H$, choose the representatives $r_x$ in
$R$, and apply \eqref{cosetformula} in $R$ to the given derivation on
$J$. This gives a derivation $\Delta:R\to T$. Each $u\in K$ fixes every
coset $x$, so $r_{ux}=r_x$ and $\Delta(u)(x)=s(r_x^{-1}u r_x)^{r_x^{-1}}=s(u)$.
Thus $\Delta|_K=\eta s$. Since $K$ acts trivially on $M$ and fixes
every coset, it acts trivially on $T$. Modulo $\eta M$ the derivation
vanishes on $K$, so it factors through $L=R/K$. By
Lemma~\ref{restriction}, it lifts to a derivation $c:L\to T$.
Let $\bar c$ be the composite of $c$ with $R\to L$. Then $\Delta-\bar c$
is a derivation into $\eta M$ whose restriction to $K$ is $\eta s$,
as desired.
\end{proof}

\begin{corollary}\label{split}
Let $L$ be a finite $p$-group, $H\le L$, and let $R$ be an extension
of $L$ by the $kL$-module $M$ inducing the given action. Assume that
$b_H(M)=[L:H]b_L(M)$. If the extension splits over $H$, it splits
over $L$. Equivalently, restriction $H^2(L,M)\to H^2(H,M)$ is
injective.
\end{corollary}
\begin{proof}
Let $C$ be a complement to $M$ in the inverse image of $H$.
The formula $d(cm)=m$ for $c\in C$ and $m\in M$ defines a derivation,
since $(cm)(c'm')=cc'(m^{c'}+m')$ in additive notation on $M$.
It restricts to the identity on $M$. Lemma~\ref{extend} extends this
restriction to a derivation $\delta:R\to M$. Its kernel meets $M$
trivially and meets every coset $gM$ in $g(-\delta(g))$. Thus
$\ker\delta$ is a complement to $M$ in $R$, as desired. The last
assertion is a reformulation, since a class in $H^2(L,M)$ is zero
precisely when the corresponding extension splits, and its restriction
to $H$ is zero precisely when the extension splits over $H$
\cite[Ch.~IV]{Brown}.
\end{proof}

\begin{lemma}\label{regular}
Let $L$ be a finite $p$-group.
\begin{enumerate}
\item[\textup{(i)}] We have $h_L(M)=0$ if and only if $M$ is a
free $kL$-module. In this case, the free rank is $\dim M^L$.
\item[\textup{(ii)}] If $h_L(M)=0$, every extension of $L$ with
kernel $M$ inducing this action splits.
\item[\textup{(iii)}] If $p$ is odd and $M$ is a nonzero free
$kL$-module, then $(g-1)^2\ne0$ on $M$ for every $1\ne g\in L$.
\end{enumerate}
\end{lemma}
\begin{proof}
Set $V_0=M^L$ and let $\rho:M\to V_0$ be a linear retraction.
The map $v\mapsto (x\mapsto\rho(v^x))$ embeds $M$ in the functions
$L\to V_0$ with action by left translation on the argument. Its
kernel has zero intersection with $V_0$, and every nonzero submodule
has a nonzero fixed point. In particular $\dim M\le |L|\dim V_0$.
If $h_L(M)=0$, Lemma~\ref{restriction} with $H=1$ gives the
reverse inequality, since $b_1(M)=-\dim M$ and $b_L(M)=-\dim V_0$.
The embedding is then onto, so $M$ is free of rank $\dim V_0$.
Conversely, a free $kL$-module of rank $r$ is isomorphic to the
functions $L\to k^r$ with action $f^g(x)=f(gx)$. By
Lemma~\ref{restriction}\textup{(i)} with $H=1$, every derivation
into this module is principal, so $h_L(M)=0$. This proves
\textup{(i)}.

For \textup{(ii)}, if $h_L(M)=0$, equality in the restriction
bound for $H=1$ gives the splitting by Corollary~\ref{split}.

For \textup{(iii)}, we may assume that $M$ is the module of functions
$L\to k^r$ with $f^g(x)=f(gx)$. Let $f$ be supported at $1$ with
$f(1)=v\ne0$, and let $1\ne g\in L$. Then $f^g$ and $f^{g^2}$ are
supported at $g^{-1}$ and $g^{-2}$, respectively. Since
$g^{-1}\notin\{1,g^{-2}\}$, we obtain that
\[
 \bigl(f^{g^2}-2f^g+f\bigr)(g^{-1})=-2v\ne0
\]
as $p$ is odd. Thus $(g-1)^2\ne0$ on $M$. This completes the proof.
\end{proof}

\begin{lemma}\label{iterate}
Let $\alpha$ be an automorphism of a group $P$ inducing the identity
on $P/E$ for an $\alpha$-invariant normal subgroup $E$. Assume that
$\alpha(e)=e\,s(e)$ for $e\in E$, where $s$ is an endomorphism of
$E$ with $s(s(e))=1$ for all $e\in E$. Writing $\alpha(g)=g e_g$, we have
\[
 g^{-1}\alpha^n(g)=e_g^n s(e_g)^{\binom n2}.
\]
If $E$ has exponent dividing an odd prime $p$, then $\alpha^p=1$,
and $\alpha$ has order $p$ when $s$ is nontrivial. If $E$ has exponent
dividing two, then $\alpha^2(g)=g s(e_g)$ and $\alpha^4=1$.

When $E$ is abelian and written additively, the hypotheses on
$\alpha|_E$ and $s$ read $\alpha|_E=1+s$ and $s^2=0$.
\end{lemma}
\begin{proof}
Let $e,f\in E$. Since $\alpha(ef)=\alpha(e)\alpha(f)$ and $s$ is an
endomorphism, we have
$$
 ef\,s(e)s(f)=e\,s(e)f\,s(f),
$$
and so $fs(e)=s(e)f$. Thus $s(E)\le Z(E)$.

We prove the formula by induction on $n$. The case $n=1$ is the
definition of $e_g$. Set $e=e_g$ and $w=e^n s(e)^{\binom n2}$, and
assume that $\alpha^n(g)=gw$. Since $s(s(e))=1$, we have $s(w)=s(e)^n$.
As $s(e)$ is central in $E$, we obtain that
$$
 \alpha^{n+1}(g)=\alpha(g)\alpha(w)=ge\,w\,s(w)
 =g\,e^{n+1}s(e)^{\binom n2+n}=g\,e^{n+1}s(e)^{\binom{n+1}2}.
$$
This completes the induction.

If $E$ has exponent dividing an odd prime $p$, then $e^p=1$ and $p$
divides $\binom p2$. Thus $\alpha^p=1$. If $s$ is nontrivial, then
$\alpha\ne1$, and so $\alpha$ has order $p$. If $E$ has exponent
dividing two, then $n=2$ gives $\alpha^2(g)=g\,s(e_g)$, and $n=4$
gives $\alpha^4(g)=g\,e_g^4s(e_g)^6=g$. This completes the proof.
\end{proof}

\section{The centralizer of a maximal normal elementary abelian subgroup}\label{sec:power}

Let $P$ be a finite $p$-group which is not elementary abelian. Let $E$
be maximal among its normal elementary abelian subgroups. Since
$E\Omega_1(Z(P))$ is normal and elementary abelian, it is equal to $E$.
Write
$$
 Z=Z(P),\quad V=\Omega_1(Z),\quad G=P/E,\quad D=C_P(E),\quad
 Y=\{y\in D:[P,y]\le E\}.
$$
Conjugation by $y$ belongs to $D_E(P)$ precisely when $y\in Y$.
In particular, $Y/E=(D/E)\cap Z(G)$. Note also that $E^G=V$.
We use right modules and take dimensions over $\F$. With
$e^g=g^{-1}eg$ and $[g,h]=g^{-1}h^{-1}gh$, a derivation into
$E$ satisfies $\delta(gh)=\delta(g)^h+\delta(h)$. The results of
this section hold for every prime $p$, except where $p$ is assumed
to be odd.

\begin{proposition}\label{power}
With the notation above, we have $Y'\le V$ and $Y^p\le Z$, and every
element of order $p$ in $Y$ lies in $E$. Moreover
$$
 [D_E(P):D_E(P)\cap\Inn(P)]
       =\frac{p^{h_G(E)}}{|Y:EZ|}.
$$
In particular, $E$ is saturated if and only if $p^{h_G(E)}=|Y:EZ|$.
If $p$ is odd, the $p$th power map on $Y$ is a homomorphism into $Z$
with kernel $E$, and saturation implies $h_G(E)\le\dim V$, with
equality precisely when $Y^p=Z$.
\end{proposition}
\begin{proof}
Let $y\in Y$ and $g\in P$. The commutator $[g,y]\in E$ commutes with
$y$, so $[g,y^p]=[g,y]^p=1$. Thus $Y^p\le Z$. Also $Y'\le E\le Z(Y)$.
Conjugation by $g$ changes each element of $Y$ by a factor in
$E\le Z(Y)$, and so it fixes every commutator of two elements of $Y$.
Hence $Y'\le E\cap Z=V$. If $y^p=1$, then $E\langle y\rangle$
is normal and elementary abelian, so $y\in E$.

It remains to count the elements of $D_E(P)$. The correspondence
$\alpha\mapsto(gE\mapsto g^{-1}\alpha(g))$ identifies $D_E(P)$ with
$\Der(G,E)$, which has order $p^{h_G(E)}|E:V|$. Conjugation by $Y$
gives exactly $|Y:Z|$ inner automorphisms in $D_E(P)$. Since
$E\cap Z=V$, we have $|Y:Z|=|Y:EZ|\,|E:V|$, which yields the formula.

Now let $p$ be odd. Since $Y'\le V$ has exponent $p$, the class-two
power formula makes the $p$th power map a homomorphism on $Y$. Its
kernel is $E$ by the first paragraph. Hence $|Y:E|=|Y^p|$ and
$|EZ:E|=|Z:V|=|Z^p|$, and so $|Y:EZ|=|Y^p:Z^p|$. If $E$ is saturated,
we obtain $p^{h_G(E)}=|Y^p:Z^p|\le |Z:Z^p|=|V|$, with equality
precisely when $Y^p=Z$. This completes the proof.
\end{proof}

\begin{proposition}\label{free}
Saturation with $Y=EZ$ is equivalent to $P=E\rtimes G$ with $E$ free
as an $\F G$-module. In this case $D=Y=E$, and the free rank is
$\dim V$. In particular, a saturated $E$ with $D=E$ always has this
form. If $p$ is odd, the condition $Y=EZ$ is equivalent to $Y^p=Z^p$.

Under these equivalent conditions, put $W=[E,P]$. Then $V<E$
and $V\le W$. Every linear map $s:E\to V$ vanishing on $W$ extends to a homomorphism
$\chi:P\to V$ trivial on a chosen complement to $E$. The map
$g\mapsto g\chi(g)$ is an automorphism inducing the identity on $P/W$
and restricting to $1+s$ on $E$.
\end{proposition}
\begin{proof}
If $E$ is saturated and $Y=EZ$, then $h_G(E)=0$ by
Proposition~\ref{power}. Lemma~\ref{regular} identifies $h_G(E)=0$
with freeness and gives the splitting and the free rank. Conversely,
if $E$ is free, it is faithful since it is nonzero, so $D=E$. Then
$Z\le D=E$ and $Y=E=EZ$, and the index formula in
Proposition~\ref{power} shows that $E$ is saturated, since
$h_G(E)=0$. The same argument shows that $Y=E=EZ$ whenever $D=E$.
If $p$ is odd, the $p$th power map on $Y$ has kernel $E$ by
Proposition~\ref{power}, and so $Y=EZ$ if and only if $Y^p=Z^p$.

Since $P$ is not elementary abelian, $E<P$ and $G\ne1$. The
regular action gives $V<E$. On each regular summand the constant
functions belong to the augmentation submodule, since $p$ divides
$|G|$. Thus $V\le W$.
Since $s$ vanishes on $W$, it is invariant under $G$ and extends to
a homomorphism $\chi:P\to V$ trivial on the chosen complement.
As $V$ is central and contained in $W$, the map $g\mapsto g\chi(g)$
is an endomorphism inducing the identity on $P/W$ and hence on
$P/E$. Its restriction $1+s$ to $E$ is invertible since $s(V)=0$,
so it is an automorphism. This completes the proof.
\end{proof}

\begin{example}\label{wreath}
Let $p$ be odd, let $P=C_p\wr C_p$, and take $E$ to be its base group.
Then $E$ is a regular module for the complement $\langle x\rangle$ of
order $p$. In particular $C_P(E)=E$, so $E$ is maximal normal
elementary abelian, and Proposition~\ref{free} shows that $E$ is
saturated. We identify $E$ with $\F[u]/(u^p)$, where $u=x-1$. Here
$uE$ is the unique invariant hyperplane and $V=u^{p-1}E$. The map
$e\mapsto e+u^{p-1}e$ extends by the identity on $\langle x\rangle$
to an automorphism in $D_{uE}(P)$. Its action on $E$ is nontrivial
and $\bigl(u^{p-1}\bigr)^2=0$. An inner automorphism acts on $E$
through $P/E$, so Lemma~\ref{regular}\textup{(iii)} shows that the
automorphism is noninner. Thus passing from $E$ to $uE$ is necessary
in this example.
\end{example}

\begin{proposition}\label{structure}
Assume that $E$ is saturated, $D>E$ and $b_G(E)\le0$. The last
condition holds automatically if $p$ is odd, since $E$ is saturated,
by Proposition~\ref{power}. Then $S=D/E$ is cyclic, $D$ is abelian,
$Y^p=Z$, and $Z$ is cyclic. The group $G$ is noncyclic. The unique
subgroup $U$ of order $p$ in $S$ is central in $G$, and one has
$$
 b_G(E)=b_U(E)=0.
$$
The inverse image $N$ of $U$ in $P$ satisfies $N\le Y$, and $a^p\in V$
for every $a\in N$. For every $x\in G\setminus S$ of order $p$, the restriction of $E$
to $\F\langle x\rangle$ is free, so all Jordan blocks of $x$ on $E$
have size $p$.
\end{proposition}
\begin{proof}
The subgroup $S$ acts trivially on $E$. Set $d(S)=\dim S/\Phi(S)$.
Then Lemma~\ref{restriction} gives
$$
 0\le (d(S)-1)\dim E=b_S(E)\le[G:S]b_G(E)\le0.
$$
Thus $d(S)=1$ and $b_G(E)=0$. In particular $S$ is cyclic. Also $D$
is abelian, as $D/E$ is cyclic and $E\le Z(D)$. Hence the $p$th power
map on $Y$ is a homomorphism into $Z$, and its kernel is $E$ by
Proposition~\ref{power}. As in the proof of Proposition~\ref{power},
we get that $|Y:EZ|=|Y^p:Z^p|$. Since $E$ is saturated and
$b_G(E)=0$, Proposition~\ref{power} gives
$|Y^p:Z^p|=p^{h_G(E)}=|V|=|Z:Z^p|$, and so $Y^p=Z$. Since $Y/E$
lies in the cyclic group $S$ and is isomorphic to $Y^p$, we see
that $Z$ is cyclic. Hence $\dim V=1$.

Assume that $G=\langle tE\rangle$, and let $m$ be the order of
the action of $t$ on $E$. Then $t^m$ centralizes both $t$ and $E$,
so $t^m\in Z$. It follows that $D=E\langle t^m\rangle=EZ$.
This gives $Y^p=Z^p<Z$, a contradiction. Thus $G$ is noncyclic.

Let $U$ be the subgroup of order $p$ in $D/E$ and $N$ its preimage
in $P$. As $U$ is characteristic in $D/E$, it is normal of order
$p$ in $G$, hence central. Thus $N\le Y$ and $b_U(E)=0$, since $U$
acts trivially on $E$. For $a\in N$, we have $a^p\in E$, and
$a^p\in Z$ by Proposition~\ref{power}, so $a^p\in V$. For any
$x\in G$ of order $p$ outside $D/E$, the group $U\times\langle x\rangle$
has $b_{U\times\langle x\rangle}(E)\le0$ by Lemma~\ref{restriction}.
Let $u$ generate $U$. The defining relations of
$U\times\langle x\rangle$ are $x^p=1$, $u^p=1$ and $[u,x]=1$. Since
$U$ acts trivially, the second relation imposes $p\,d(u)=0$ on a
derivation $d$, which always holds, and the third says that $d(u)$ is
fixed by $x$. If $w=d(u)$ and $d_0=d|_{\langle x\rangle}$, then
$d(u^ix^j)=i w+d_0(x^j)$. Thus a derivation on $U\times\langle x\rangle$ amounts
to a derivation on $\langle x\rangle$ and an arbitrary vector $d(u)$
fixed by $x$. Hence
\[
 b_{U\times\langle x\rangle}(E)
 =\dim\Der(\langle x\rangle,E)+\dim E^{\langle x\rangle}-\dim E
 =h_{\langle x\rangle}(E).
\]
This is nonnegative, so it vanishes. Lemma~\ref{regular} makes $E$
free over $\F\langle x\rangle$, as required.
\end{proof}

\begin{corollary}\label{noncyclic-centre}
Let $p$ be odd, and assume that $Z(P)$ is noncyclic. Then $E$ is
saturated if and only if $P$ splits over $E$ and $E$ is free as an
$\F[P/E]$-module.
\end{corollary}
\begin{proof}
If $E$ is saturated, Proposition~\ref{structure} excludes $D>E$.
Proposition~\ref{free} now gives the assertion in both directions.
\end{proof}

\section{Lifting transvections}\label{sec:lift}

\begin{lemma}\label{inner-lift}
Let $p$ be odd. Assume that $E$ is saturated and $\alpha$ is an inner
automorphism inducing the identity on $P/E$. If $\alpha|_E=1+s$ for a
linear map $s:E\to E$ with $s^2=0$, then $\alpha$ fixes $E$ pointwise.
\end{lemma}
\begin{proof}
Let $y\in P$ be such that $\alpha$ is conjugation by $y$. Then
$yE\in Z(G)$, and Lemma~\ref{iterate} gives $y^p\in Z$. First assume
that $D=E$. Then $E$ is free over $\F G$ by Proposition~\ref{free}.
Since $(yE-1)^2=0$ on $E$, Lemma~\ref{regular}\textup{(iii)} gives
$yE=1$. Thus $y\in E$, and so $\alpha$ fixes $E$ pointwise.

Now assume that $D>E$. Since $p$ is odd, we have $b_G(E)\le0$ by
Proposition~\ref{power}. By Proposition~\ref{structure}, we have
$Y^p=Z$, and so there is $j\in Y$ with $j^p=y^p$. Since $yE$ is
central in $G$, the element $q=yj^{-1}E$ has order dividing $p$.
It acts as $y$ on $E$ because $j$ centralizes $E$. Assume that this
action is nontrivial. Then $q$ has order $p$ and lies outside $S$, so
$E$ is free over $\F\langle q\rangle$ by Proposition~\ref{structure}.
Lemma~\ref{regular}\textup{(iii)} then gives $(q-1)^2\ne0$ on $E$,
which is a contradiction as $q$ acts as $1+s$ with $s^2=0$. Thus $y$
acts trivially on $E$, as desired.
\end{proof}

\begin{lemma}\label{lift-transvection}
Assume that $E$ is saturated and $b_G(E)\le0$. For every
$G$-equivariant linear map $s:E\to V$ vanishing on $V$, there is an
automorphism of $P$ inducing the identity on $P/E$ and restricting to
$1+s$ on $E$.
\end{lemma}
\begin{proof}
The map $s$ vanishes on $[E,P]$ since $V$ is fixed by $G$.
If $D=E$, the lift is given by Proposition~\ref{free}.

Now assume that $D>E$, and let $U$ and $N$ be as in
Proposition~\ref{structure}. Write $N=E\langle a\rangle$. Then
$N\le Y$ and $a^p\in V$.
Since $s(a^p)=0$, the assignments $d|_E=s$ and $d(a)=0$ define a
derivation $N\to E$, where $N$ acts trivially. As $b_U(E)=b_G(E)=0$,
Lemma~\ref{extend} extends $s$ to a derivation $\delta:P\to E$.
The endomorphism $g\mapsto g\delta(g)$ induces the identity on
$P/E$ and the invertible map $1+s$ on $E$, so it is an automorphism,
as desired.
\end{proof}

The next lemma collects further consequences of
Proposition~\ref{structure} that we use for every prime. Part
\textup{(ii)} is a transfer argument.

\begin{lemma}\label{transfer}
Assume that $E$ is saturated, $D>E$ and $b_G(E)\le0$. Let $U$ and $N$
be as in Proposition~\ref{structure}, and put $W=[E,P]$.
\begin{enumerate}
\item[\textup{(i)}] We have $V\le W$.
\item[\textup{(ii)}] Every derivation $G\to E$ takes values in $W$.
In particular, $[P,N]\le W$.
\item[\textup{(iii)}] Let $\alpha:g\mapsto g\delta(g)$ be an
automorphism of $P$ inducing the identity on $P/E$. If $\alpha$ fixes
$W$ pointwise, then $\delta(N)\le V$.
\item[\textup{(iv)}] Let $\alpha:g\mapsto g\delta(g)$ be an
automorphism of $P$ inducing the identity on $P/E$ and fixing $W$
pointwise, and let $U_2=U\times\langle x\rangle$ be a normal subgroup
of $G$ with $x$ of order $p$. Then, after composing $\alpha$ with conjugation by a
suitable element of $E$, we have $\delta(g)(x-1)\in V$ for all
$g\in P$. This does not change $\alpha|_E$.
\end{enumerate}
\end{lemma}
\begin{proof}
By Proposition~\ref{structure}, $G$ is noncyclic and $S$ is cyclic.
Thus $G\ne S$, and so $W\ne0$. The nonzero submodule $W$ contains a
nonzero fixed vector. As $E^G=V$ is a line by
Proposition~\ref{structure}, we get that $V\le W$. This proves (i).

We use transfer to show that every derivation $c:G\to E$ takes
values in $W$. Since $b_U(E)=b_G(E)=0$, Lemma~\ref{restriction}
gives $c=\varepsilon\Delta$ with $\Delta$ as in \eqref{cosetformula}
for some $d:U\to E$, up to a principal derivation. Since $E/W$ has
trivial action, the projection $\pi:E\to E/W$ kills principal
derivations and gives
\[
 \pi c(g)=\sum_{\xi\in G/U}(\pi d)(r_{g\xi}^{-1}g r_\xi).
\]
The map $\pi d$ is a homomorphism. Since $U$ is normal, every orbit of
$\langle g\rangle$ on $G/U$ has length $\ell$, the order of $gU$.
Along each orbit, the arguments multiply in reverse orbit order to
$r_\xi^{-1}g^\ell r_\xi=g^\ell$, since $g^\ell\in U\le Z(G)$. As
$\pi d$ is a homomorphism into an abelian group, we obtain that
$$
 \pi c(g)=\tfrac{[G:U]}{\ell}\,(\pi d)(g^\ell)=(\pi d)\bigl(g^{[G:U]}\bigr).
$$
Since $G$ is a noncyclic $p$-group, its exponent divides
$|G|/p=[G:U]$. Hence $\pi c=0$. Applying this to the derivation
$gE\mapsto[g,y]$ for each $y\in N\le Y$, we obtain $[P,N]\le W$.

For (iii), let $y\in N$ and $g\in P$. Since $[y,g]\in W$, the
automorphism $\alpha$ fixes $[y,g]$. Applying $\alpha$ to
$y^g=y[y,g]$, we obtain that $\alpha(y^g)=y^g\delta(y)$. On the
other hand, $\alpha(y^g)=\alpha(y)^{\alpha(g)}=y^g\delta(y)^g$,
since $N$ centralizes $E$. Thus $\delta(y)$ is fixed by $P$, and so
$\delta(y)\in V$.

For (iv), note that $x\notin S$, since $U_2$ is not cyclic. By
Proposition~\ref{structure}, $E$ is free over $\F\langle x\rangle$. By
Lemma~\ref{regular}, the extension of $\langle x\rangle$ by $E$ inside
$P$ splits, and every derivation $\langle x\rangle\to E$ is principal.
Choose a lift $t$ of $x$ of order $p$. Then $\delta(t)=e^t-e$ for some
$e\in E$. Let $c_e(g)=e^{-1}ge$. Replacing $\alpha$ by
$c_e\circ\alpha$ subtracts the
principal derivation $g\mapsto e^g-e$ from $\delta$. This does not
change $\alpha|_E$, so the new automorphism still fixes $W$ pointwise,
and we may assume that $\delta(t)=0$. The $p$-group $G$ acts trivially
on the normal section $U_2/U$ of order $p$, so for each $g\in P$ there
is $b\in N$ with $t^g=tb$. Since $\alpha(t)=t$ and $\delta(b)\in V$ by
\textup{(iii)}, we obtain that $\alpha(t^g)=tb\delta(b)=t^g\delta(b)$.
On the other hand, $\alpha(t^g)=\alpha(t)^{\alpha(g)}=(t^g)^{\delta(g)}$.
Since $b$ centralizes $E$, conjugation by $t^g$ acts as $x$ on $E$, and
so $(t^g)^{\delta(g)}=t^g\cdot\delta(g)(1-x)$ in additive notation on
$E$. Hence $\delta(g)(x-1)\in V$. This completes the proof.
\end{proof}

\section{Fixing a hyperplane and its quotient}\label{sec:flag}

In this section, $p$ is odd. To control the action of the lift on
$P/F$, we use a normal subgroup of order $p^2$ in $P/E$.

\begin{lemma}\label{plane}
Let $G$ be a noncyclic finite $p$-group for odd $p$, and let
$U\le Z(G)$ have order $p$. There is a normal elementary abelian
subgroup $U_2$ of order $p^2$ containing $U$.
\end{lemma}
\begin{proof}
By \cite[Theorem~5.4.10(i)]{Gor}, the noncyclic group $G$ has a
noncyclic abelian normal subgroup $A$. Then $\Omega_1(A)$ is a normal
elementary abelian subgroup of order at least $p^2$, and so it contains
a normal subgroup $A_2$ of $G$ of order $p^2$. If $U\le A_2$, we put
$U_2=A_2$. Otherwise, we choose a subgroup $L$ of order $p$ in
$A_2\cap Z(G)$, which is nontrivial as $A_2$ is normal, and we put
$U_2=U\times L$. In both cases $U_2$ is as desired.
\end{proof}

\begin{proposition}\label{displacement}
Assume that $E$ is saturated, and put $W=[E,P]$. Then $V\le W$.
For every linear map $s:E\to V$ vanishing on $W$, there is an
automorphism of $P$ inducing the identity on $P/E$ and restricting
to $1+s$ on $E$. Every such lift induces the identity on $P/W$.
\end{proposition}
\begin{proof}
Two lifts of $1+s$ inducing the identity on $P/E$ give derivations
$P\to E$ that agree on $E$. Their difference therefore factors
through a derivation $G\to E$. If every such derivation takes values
in $W$, one lift inducing the identity on $P/W$ suffices to prove
the assertion for all lifts.

If $D=E$, Proposition~\ref{free} gives $V\le W$ and a lift inducing
the identity on $P/W$. The module $E$ is free, so every derivation
$G\to E$ is principal by Lemma~\ref{regular} and takes values in $W$.

Now assume that $D>E$. Since $p$ is odd, we have $b_G(E)\le0$ by
Proposition~\ref{power}, so Proposition~\ref{structure} and
Lemma~\ref{transfer} apply. Thus $V\le W$, and every derivation
$G\to E$ takes values in $W$. Since $s$ vanishes on $W$ and $V\le W$,
it is $G$-equivariant and vanishes on $V$. Let $\alpha(g)=g\delta(g)$
be a lift given by Lemma~\ref{lift-transvection}. It fixes $W$
pointwise, as $s(W)=0$.

Let $U_2$ be as in Lemma~\ref{plane}, and choose $x\in U_2\setminus U$.
Then $x$ has order $p$ and $U_2=U\times\langle x\rangle$. By
Lemma~\ref{transfer}\textup{(iv)}, we may assume that
$\delta(g)(x-1)\in V$ for all $g\in P$. Since $V$ is fixed by $x$, we
get that $\delta(g)(x-1)^2=0$. By Proposition~\ref{structure}, $E$ is
free over $\F\langle x\rangle$. In the regular module $\F[u]/(u^p)$
with $u=x-1$, the kernel of $u^2$ is generated by $u^{p-2}$, and so it
lies in the augmentation submodule as $p>2$. Hence
$\delta(g)\in E(x-1)\le W$. Thus this lift induces the identity on
$P/W$, and the first paragraph completes the proof.
\end{proof}

\begin{proof}[Proof of Theorem~\ref{flag}]
Put $W=[E,P]$. Proposition~\ref{displacement} gives $V\le W$,
and $W<E$ because a simple quotient of the nonzero $G$-module
$E$ is trivial. Let $F$ be a $P$-invariant hyperplane of $E$.
The normal subgroup $E/F$ of order $p$ is central in $P/F$,
so $V\le W\le F$. Let $s:E\to V$ be a linear map with kernel $F$ and
one-dimensional image. By Proposition~\ref{displacement}, $1+s$
has a lift inducing the identity on $P/E$ and on $P/W$. This lift
fixes $F$ and belongs to $D_F(P)$. Since $s\ne0$ and $s^2=0$, it
is noninner by Lemma~\ref{inner-lift} and has order $p$ by
Lemma~\ref{iterate}. This completes the proof.
\end{proof}

\begin{proof}[Proof of Theorem~\ref{berkovich}]
Let $E$ be maximal among the normal elementary abelian subgroups of
$P$. First assume that $E$ is unsaturated, and let $\alpha\in D_E(P)$
be noninner. Then $\alpha$ has order $p$, and it is quadratic since
every element of $D_E(P)$ is. Now assume that $E$ is saturated. As in the
proof of Theorem~\ref{flag}, we have $[E,P]<E$. Since $P$ acts
trivially on $E/[E,P]$, every hyperplane $F$ of $E$ containing
$[E,P]$ is $P$-invariant. Theorem~\ref{flag} gives a noninner
$\alpha\in D_F(P)$ of order $p$, and $\alpha$ is quadratic as an
element of $D_F(P)$. This completes the proof.
\end{proof}

\begin{corollary}\label{outer-interval}
Let $p$ be odd. Assume that $E$ is saturated and put $W=[E,P]$. For every $W\le F\le E$,
the map $\alpha\mapsto(eF\mapsto[e,\alpha])$ induces an isomorphism
$$
 D_F(P)/(D_F(P)\cap\Inn(P))\ \cong\ \Hom_{\F}(E/F,V).
$$
For $F=W$ the image of $D_W(P)$ in $\operatorname{Out}(P)$ is elementary abelian
of rank $\dim(E/W)\dim V$.
\end{corollary}
\begin{proof}
Since $P$ acts trivially on $E/W$, the subgroup $F$ is normal.
Proposition~\ref{displacement} gives $V\le W\le F$. For
$\alpha\in D_F(P)$, its displacement on $E$ induces a $P$-map
$E/F\to F$ with image in $F^P=V$. These maps add under composition.
An element of the kernel belongs to $D_E(P)$, so it is inner by
saturation. Conversely, an inner element of $D_F(P)$ acts on $E$
as $1+s$ with $s^2=0$, because $s(E)\le F$ and $s(F)=0$.
By Lemma~\ref{inner-lift} it fixes $E$ pointwise.

For any $s:E/F\to V$, the composite $E\to E/F\to V$ vanishes on
$W$. Proposition~\ref{displacement} gives a lift that fixes $F$
and induces the identity on $P/W$, hence belongs to $D_F(P)$
and maps to $s$. This proves surjectivity. Taking $F=W$, we obtain
the stated rank. This completes the proof.
\end{proof}

We next relate this interval of subgroups to the construction in
\cite[Sections~2.2--2.3]{Kizmaz}. Recall the definition, which
applies also to $p=2$. For an abelian normal subgroup $B$ of $P$, put
$$
 Y(B)=\{y\in C_P(B):[P,y]\le B\},\qquad
 C(B)=\{y\in Y(B):y^p\in B\}.
$$
We have $Y(B)/B\le Z(P/B)$, and $C(B)$ is the inverse image of
$\Omega_1(Y(B)/B)$ in $Y(B)$. In \cite[Section~2.2]{Kizmaz}, the
subgroup $\mathcal C(B)$ is defined as the subgroup generated by the
elements $x\in C_P(B)$ with $x^p\in B$ and $[P,x]\le B$. Since
$Y(B)/B$ is abelian, these elements form a subgroup, and so
$\mathcal C(B)=C(B)$. Given a normal elementary abelian subgroup $F$,
start with $M_0=Z(Y(F))$ and put $M_{i+1}=Z(C(M_i))$. Since each
abelian $M_i$ is contained in $Z(C(M_i))$, the series stabilizes at a
subgroup $M$. We write $K_P(F)=C(M)$. In \cite[Section~2.3]{Kizmaz},
for an abelian normal subgroup $U$ of $P$, one puts
$\mathcal M(P,U)=Z(\mathcal C(U))$, iterates
$\mathcal M_k=\mathcal M(P,\mathcal M_{k-1})$ starting from
$\mathcal M_1=\mathcal M(P,U)$, and defines
$\mathcal K(P,U)=\mathcal C(\mathcal M_\infty)$ for the terminal term
$\mathcal M_\infty$. Thus $K_P(F)=\mathcal K(P,Z(Y(F)))$.

For odd $p$ and saturated $E$, all the subgroups $F$ in
Corollary~\ref{outer-interval} give the same group. In fact,
$$
 K_P(F)=C_P(E)\qquad\text{for every }[E,P]\le F\le E .
$$
To see this, note that $E\le Y(F)$. Conjugation by $y\in Y(F)$ belongs to $D_F(P)$
and its restriction to $E$ is $1+s$ with $s^2=0$. It fixes $E$
pointwise by Lemma~\ref{inner-lift}, so $Y(F)\le D=C_P(E)$.
Thus $E\le M_0$. As $D$ is abelian, at each stage
$E\le M_i\le C(M_i)\le D$, so $Z(C(M_i))=C(M_i)$.
If $M<D$, the nontrivial normal subgroup $D/M$ of $P/M$ contains a
central subgroup of order $p$. A lift of a generator centralizes $M$
and belongs to $C(M)\setminus M$, contrary to stability.

Propositions~\ref{free} and~\ref{structure} leave the following
question for odd primes. Can a maximal normal elementary abelian
subgroup $E$ be saturated with $C_P(E)>E$? By the preceding
paragraph, this is equivalent to asking whether a saturated $E$ can
satisfy $K_P(E)>E$. Any such example must have cyclic center.
Using GAP and the Small Groups Library \cite{GAP,SmallGrp}, we checked all pairs $(P,E)$
with $P$ nonabelian of order $3^4$, $3^5$, $3^6$, $5^4$, $5^5$,
$7^4$ or $7^5$ and $E$ maximal among the normal elementary abelian
subgroups of $P$. Up to isomorphism, the only saturated pair
consists of $C_3\wr C_3$ and its base group $E$, for which
$C_P(E)=E$.

\section{Two-groups: a quaternion subgroup and the split case}\label{sec:two}

In Sections~\ref{sec:two}--\ref{sec:Kp}, we have $p=2$. Let $P$ be a
finite nonabelian $2$-group, and let $E$, $Z$, $V$, $G$, $D$ and $Y$
be as in Section~\ref{sec:power}. By Lemma~\ref{iterate} with $s=0$,
every element of $D_E(P)$ has order at most two, and it is quadratic.
Thus if $E$ is unsaturated, $D_E(P)$ contains a quadratic noninner
involution. We may therefore assume until the end of
Section~\ref{sec:cyclic} that $E$ is saturated. Put
$$
 n=\dim E,\qquad z=\dim V .
$$
By Proposition~\ref{power}, we have $2^{h_G(E)}=|Y:EZ|$, and so
$b_G(E)=\log_2|Y:EZ|-z$.

We also use automorphisms fixing an abelian normal subgroup $B$ of $P$
which need not be elementary abelian. Let $D_B(P)$ be the group of
automorphisms of $P$ fixing both $B$ and $P/B$ pointwise. Each
$\alpha\in D_B(P)$ has the form $g\mapsto gf(gB)$ for a unique map
$f:P/B\to B$, and the homomorphism condition for $\alpha$ is exactly
the identity $f(gh)=f(g)^hf(h)$, written multiplicatively in $B$.
We call such maps derivations as well. An inner automorphism belongs
to $D_B(P)$ only if it is conjugation by an element of $C_P(B)$.

\begin{proposition}\label{large}
Assume that $b_G(E)>0$. Then $|Z|=2$, $h_G(E)=2$, and $D=Y=EJ$ for a
subgroup $J\cong Q_8$ with $E\cap J=Z$. Moreover, $E$ is a free
$\mathbb F_2[P/Y]$-module of rank one.
\end{proposition}
\begin{proof}
Set $T=Y/E$, $Q=P/Y$, $A=Z(Y)$ and $h=h_G(E)$. Since $E$ is
saturated, every derivation $G\to E$ is induced by conjugation by an
element of $Y$, which is determined modulo $Z$. The principal
derivations are induced by $E$. Since $T$ acts trivially on $E$, such
a derivation vanishes on $T$ precisely when the conjugating element
centralizes $Y$, that is, lies in $A$. The derivations vanishing on
$T$ are exactly those inflated from $Q$. Consequently
$$
 \nu:=h_Q(E)=\log_2|A:EZ|.
$$
By Proposition~\ref{power}, squaring on the abelian group $A$ is a
homomorphism into $Z$ with kernel $E$. Hence
$2^\nu=|A^2:Z^2|\le|Z:Z^2|=2^z$. Assume that $\nu=z$. Then $A^2=Z$.
For each $y\in Y$ we can choose $c\in A$ with $c^2=y^2$. Since $c$
centralizes $y$, we get $(yc^{-1})^2=1$, and hence
$yc^{-1}\in E\le A$. Thus $Y=A$, and so $h=\nu=z$, which contradicts
$b_G(E)>0$. Hence $\nu<z$.

Since $T$ acts trivially on $E$, we have $b_T(E)=(d(T)-1)n$, where
$d(T)$ is the minimal number of generators of $T$. Moreover
$E^Q=E^G=V$, and so $b_Q(E)=\nu-z$. Applying
Lemma~\ref{restriction} to $1\le Q$ and to $T\le G$, we obtain that
$$
 -n\le|Q|(\nu-z),\qquad (d(T)-1)n\le|Q|(h-z).
$$
The group $Y/EZ$ is elementary abelian of rank $h$, since $Y^2\le Z$
and $Y'\le V$, and it is a quotient of $T$. Hence $d(T)\ge h$.
Moreover, $h\ge2$ as $h>z\ge1$. Consequently,
$$
 (h-1)|Q|(z-\nu)\le(h-1)n\le(d(T)-1)n\le|Q|(h-z),
$$
that is, $(h-1)(z-\nu)\le h-z$. As $z-\nu\ge1$, this gives
$h-1\le h-z$, and so $z=1$ and $\nu=0$. Thus $A=EZ$, and $E$ is a free
$\mathbb F_2Q$-module of rank $\dim E^Q=1$ by Lemma~\ref{regular}. A
nonzero free module is faithful, so the kernel $D/Y$ of the action of
$Q$ on $E$ is trivial, that is, $D=Y$.

Since $z=1$, the group $Z$ is cyclic. Write $V=\langle v\rangle$. Let
$y\in Y\setminus EZ$. If $y^2=c^2$ for some $c\in Z$, then
$(yc^{-1})^2=1$, and so $yc^{-1}\in E$, which is not the case. Since
$|Z:Z^2|=2$, it follows that $y^2$ generates $Z$. Now let $x,y\in Y$
be independent modulo $EZ$. Since $Y'\le V\le Z(Y)$, we have
$(xy)^2=x^2y^2[y,x]$. Assume that $|Z|\ge4$. Then $V\le Z^2$, and
$x^2y^2\in Z^2$ as $x^2$ and $y^2$ both generate $Z$. Thus
$(xy)^2\in Z^2$, which is not possible as $xy\notin EZ$. Hence $Z=V$,
and $x^2=y^2=(xy)^2=v$. It follows that $[y,x]=v$ whenever $x$ and $y$
are independent modulo $EZ$. If $Y/EZ$ contained the images of three
independent elements $x_1,x_2,x_3$, we would get that
$$
 [x_1,x_2x_3]=[x_1,x_3][x_1,x_2]=v^2=1,
$$
although $x_1$ and $x_2x_3$ are independent modulo $EZ$. Hence
$h=2$. Note that $EZ=E$. Let $y,t\in Y$ be lifts of a basis of $Y/E$.
Then $y^2=t^2=[y,t]=v$, and so $J=\langle y,t\rangle$ is a quaternion
group of order eight. Since $J/Z$ maps isomorphically onto $Y/E$, we
have $E\cap J=Z$ and $Y=EJ$, as desired.
\end{proof}

\begin{lemma}\label{invert}
Let $B$ be an abelian normal subgroup of $P$ of exponent at most four,
and let $t\in P$ be an element which inverts $B$ and satisfies
$[P,t]\le B$. Let $c_t$ denote conjugation by $t$. Then $c_t\gamma$ is
a cubic automorphism of order at most two for every
$\gamma\in D_B(P)$. Moreover, $c_t\gamma$ is noninner whenever
$\gamma$ is noninner.
\end{lemma}
\begin{proof}
Put $\beta=c_t\gamma$. Since $c_t$ inverts $B$ and $\gamma$ fixes $B$
pointwise, $\beta$ inverts $B$. Since $[P,t]\le B$, both $c_t$ and
$\gamma$ induce the identity on $P/B$, and hence so does $\beta$. Let
$g\in P$, and write $\beta(g)=gb$ with $b\in B$. Then
$$
 \beta^2(g)=\beta(g)\beta(b)=gbb^{-1}=g.
$$
Moreover, $[P,\beta]\le B$, $[B,\beta]=B^2$ and $[B^2,\beta]=B^4=1$,
as $B$ has exponent at most four. Thus $\beta$ is cubic. The last
assertion holds since $c_t$ is inner. This completes the proof.
\end{proof}

\begin{lemma}\label{order-eight}
Every nonabelian group of order eight has a cubic noninner involution.
\end{lemma}
\begin{proof}
Let $P$ be nonabelian of order eight, and let $B=\langle b\rangle$ be
a cyclic subgroup of order four. Then $B$ is normal of index two. Let
$t\in P\setminus B$. Since $P$ is nonabelian, $t$ acts nontrivially on
$B$, and so it inverts $B$. Also $[P,t]\le B$, as $P/B$ has order two.
The map $f:P/B\to B$ with $f(B)=1$ and $f(tB)=b$ is a derivation,
since $f(tB)^tf(tB)=b^{-1}b=1$. Let $\gamma\in D_B(P)$ be the
corresponding automorphism. Then $\gamma(t)=tb$. An inner
automorphism in $D_B(P)$ is conjugation by an element $c$ of
$C_P(B)=B$, and it sends $t$ to $t[t,c]=tc^2$. Since $b\notin B^2$,
the automorphism $\gamma$ is noninner. By Lemma~\ref{invert},
$c_t\gamma$ is a cubic noninner involution, as desired.
\end{proof}

\begin{proposition}\label{invertcase}
If $b_G(E)>0$, then $P$ has a cubic noninner involution.
\end{proposition}
\begin{proof}
Let $J=\langle y,t\rangle$ be as in Proposition~\ref{large}, and set
$Q=P/Y$ and $B=E\langle y\rangle$. Since $y\in Y$, the subgroup $B$ is
normal in $P$. It is abelian of exponent four as $y^2\ne1$. Since $t$
centralizes $E$ and $y^t=y[y,t]=y^{-1}$, the element $t$ inverts $B$.
Moreover, $[P,t]\le E\le B$ as $t\in Y$.

First assume that $Q\ne1$. We identify $B/E$ with $\mathbb F_2$ via
$yE$. Choose a nonzero character $\chi:Q\to B/E$, and inflate it to
$L=P/B$. Consider the short exact sequence $0\to E\to B\to B/E\to0$ of
abelian groups with $L$-action. Here $L$ acts trivially on $B/E$. The
connecting map sends $\chi$ to a class in $H^2(L,E)$, and this class
restricts to zero on $Y/B$ since $\chi$ vanishes there. We claim that
restriction $H^2(L,E)\to H^2(Y/B,E)$ is injective. Note that $Y/B$ acts
trivially on $E$, and $E\cong\mathbb F_2[L/(Y/B)]$ by
Proposition~\ref{large}. Hence restriction to $Y/B$ followed by
evaluation at the identity coset is the isomorphism of Shapiro's
lemma \cite[III.5.9 and III.8, Exercise~2]{Brown}. Thus the connecting map sends $\chi$ to
zero, and so $\chi$ lifts to a derivation $f:L\to B$. Let
$\gamma\in D_B(P)$ be the corresponding automorphism. Assume that
$\gamma$ is conjugation by some $w\in P$. Then $w\in C_P(B)\le D=Y$,
and so $[g,w]\in E$ for all $g\in P$. This is not possible, as $f$
lifts the nonzero character $\chi$. Thus $\gamma$ is noninner, and
Lemma~\ref{invert} shows that $c_t\gamma$ is a cubic noninner
involution.

Now assume that $Q=1$. Then $\dim E=1$, so $E=Z$ and $P=Y=J\cong Q_8$.
Lemma~\ref{order-eight} completes the proof.
\end{proof}

\begin{proposition}\label{split-two}
Assume that $D=E$. Then $P$ has a quadratic noninner involution unless
$P\cong D_8$. In all cases, $P$ has a cubic noninner involution.
\end{proposition}
\begin{proof}
By Proposition~\ref{free}, we have $P=E\rtimes G$, the module $E$ is
free of rank $z$, and $V\le W=[E,P]$. Let $s:E\to V$ be a linear map
of rank one vanishing on $W$. Proposition~\ref{free} gives an
automorphism $\alpha:g\mapsto g\chi(g)$ restricting to $1+s$ on $E$,
where $\chi:P\to V$ is a homomorphism. Since $V\le W$, we have
$\chi(V)=1$, and so $\alpha^2=1$. Its displacement lies in $V$, which it fixes pointwise,
so $\alpha$ is quadratic.

Assume that $\alpha$ is conjugation by some $w\in P$. Then
$w\notin E$ as $s\ne0$. Let $k>1$ be the order of $wE$ in $G$. Then
$wE$ permutes each regular basis of $E$ in cycles of length $k$, so
$wE-1$ has rank $n(1-1/k)\ge n/2$ on $E$. Since $s$ has rank one, we
get that $n\le2$. On the other hand, $n\ge|G|\ge2$ by freeness. Thus
$n=|G|=2$, and so $P\cong C_2\wr C_2\cong D_8$. In this case,
Lemma~\ref{order-eight} gives a cubic noninner involution. This
completes the proof.
\end{proof}

\begin{example}\label{d8}
Write $D_8=\langle r,t\mid r^4=t^2=1,\ r^t=r^{-1}\rangle$. Every
automorphism of $D_8$ sends $r$ to $r^\epsilon$ and $t$ to $r^jt$ for
some $\epsilon=\pm1$ and some $j$. The four automorphisms with $j$
even are precisely the inner ones. Let $j$ be odd. Then the square of
the automorphism fixes $t$ if and only if $\epsilon=-1$. Thus every
noninner involution inverts $\langle r\rangle$ and sends $t$ to
$r^jt$, where $r^j$ generates $\langle r\rangle$. Its double
commutator is $\langle r^2\rangle$, and its triple commutator is
trivial.
\end{example}

Assume that $Z$ is noncyclic. Then $b_G(E)\le0$ by
Proposition~\ref{large}. If $D>E$, then $Z$ is cyclic by
Proposition~\ref{structure}. Hence $D=E$, and
Proposition~\ref{split-two} gives a quadratic noninner involution,
since $D_8$ has cyclic center. Together with the unsaturated case,
this proves the second assertion of Theorem~\ref{two}. The same
argument shows that Corollary~\ref{noncyclic-centre} also holds for
$p=2$.

\section{Two-groups with cyclic center}\label{sec:cyclic}

In this section, we assume that $b_G(E)\le0$ and $D>E$. By
Proposition~\ref{structure}, the group $S=D/E$ is cyclic, $D$ is
abelian, $Y^2=Z$, the center $Z$ is cyclic, the group $G$ is
noncyclic, and $b_G(E)=b_U(E)=0$, where $U=\langle u\rangle$ is the
subgroup of order two in $S$. Let $N$ be the inverse image of $U$ in
$P$, and write $V=\langle v\rangle$ and $W=[E,P]$. By
Lemma~\ref{transfer}, we have $V\le W$, every derivation $G\to E$
takes values in $W$, and $[P,N]\le W$.

\begin{lemma}\label{local}
With the notation above, the following hold.
\begin{enumerate}
\item[\textup{(i)}] Restriction $H^2(G,E)\to H^2(U,E)$ is injective.
Every element of $N\setminus E$ has square $v$.
\item[\textup{(ii)}] We have $n\ge3$. Every involution
$x\in G\setminus S$ lifts to an involution of $P$.
\end{enumerate}
\end{lemma}
\begin{proof}
Since $b_U(E)=b_G(E)=0$, Corollary~\ref{split} shows that
restriction $H^2(G,E)\to H^2(U,E)$ is injective. Let $b\in N\setminus E$.
Then $b\in Y$, so $b^2\in E\cap Z=V$ by Proposition~\ref{power}. Since
$b\notin E$, it is not an involution, and hence $b^2=v$. This proves
(i).

Since $G$ is noncyclic and $S$ is cyclic, we have $D<P$, which excludes
$n=1$. Assume that $n=2$. Then the nontrivial $2$-group $P/D$ acts
faithfully on $E$, so it has order two and is generated by a
transvection. Hence $E$ is free over $\mathbb F_2[P/D]$. Let
$y\in Y\setminus EZ$, which exists since
$|Y:EZ|=2^{h_G(E)}=2^z=2$. Since $D$ is abelian, the derivation
$gE\mapsto[g,y]$ vanishes on $D/E$, and so it is inflated from $P/D$.
By Lemma~\ref{regular}, it is principal. Thus $y$ and some element of
$E$ induce the same inner automorphism, and so $y\in EZ$, which is a
contradiction. Hence $n\ge3$.

Let $x\in G\setminus S$ be an involution. By
Proposition~\ref{structure}, $E$ is free over
$\mathbb F_2\langle x\rangle$, and so the extension of $\langle x\rangle$
by $E$ inside $P$ splits by Lemma~\ref{regular}. A complement is
generated by an involution lifting $x$. This completes the proof.
\end{proof}

\begin{lemma}\label{outer}
Let $\alpha$ be an automorphism of $P$ of order two such that
$\alpha|_E=1+s$ for a nonzero linear map $s:E\to V$. Then $\alpha$ is
noninner.
\end{lemma}
\begin{proof}
Assume that $\alpha$ is conjugation by $w\in P$. Then $w^2\in Z$ and
$w\notin D$. Since $Y^2=Z$, there is $y\in Y$ with $y^2=w^2$. As $yE$
is central in $G$, we get that $(wy^{-1})^2\in w^2y^{-2}E=E$. Thus
$x=wy^{-1}E$ is an involution of $G$ which acts on $E$ as $w$ does,
and $x\notin S$ as $w\notin D$. By Proposition~\ref{structure}, $E$ is
free over $\mathbb F_2\langle x\rangle$, and so $x-1$ has rank $n/2$ on
$E$. However, $s$ has rank one and $n\ge3$ by Lemma~\ref{local}, which
is a contradiction.
\end{proof}

\begin{lemma}\label{supplement}
Assume that $E\not\le\Phi(P)$. Then $P$ has a quadratic noninner
involution $\alpha$ with $[P,\alpha]\le V$. This applies in particular
if $P=HE$ for a proper subgroup $H$ of $P$.
\end{lemma}
\begin{proof}
By Lemma~\ref{local}, $v$ is the square of an element of $N$, and so
$V\le\Phi(P)$. Let $P_1$ be a maximal subgroup of $P$ not containing
$E$, and let $\chi:P\to V$ be the homomorphism with kernel $P_1$. Since
$\chi(P)\le V$ is central, the map $\alpha:g\mapsto g\chi(g)$ is an
endomorphism of $P$. As $V\le\Phi(P)\le P_1$, we have $\chi(V)=1$, and
so $\alpha^2=1$. Thus $\alpha$ is an involution with
$[P,\alpha]\le V$. It fixes $V$ pointwise, and so it is quadratic. Its
restriction to $E$ is $1+s$, where $s=\chi|_E$ is nonzero as
$E\not\le P_1$. By Lemma~\ref{outer}, $\alpha$ is noninner. Finally,
if $P=HE$ for a proper subgroup $H$, then a maximal subgroup containing
$H$ does not contain $E$, and so $E\not\le\Phi(P)$. This completes the
proof.
\end{proof}

\begin{proposition}\label{four}
If $G$ has a normal four-subgroup, then $P$ has a quadratic noninner
involution.
\end{proposition}
\begin{proof}
Every normal four-subgroup of $G$ meets $Z(G)$ nontrivially. If such a
subgroup does not contain $U$, then a central involution in it
generates, together with $U$, another normal four-subgroup. We may
therefore choose a normal four-subgroup $U_2=U\times\langle x\rangle$
of $G$. Then $x$ is an involution outside $S$, since $u$ is the only
involution of the cyclic group $S$. By Proposition~\ref{structure},
$E$ is free over $\mathbb F_2\langle x\rangle$. Put $\sigma_x=x+1$.
Then $\ker\sigma_x=E\sigma_x$, and this subspace has dimension $n/2$.
Each conjugate of $x$ lies in $xU$, and $U$ acts trivially on $E$.
Thus $\sigma_x$ commutes with the action of $G$ on $E$.

Put $E_0=\{e\in E:e\sigma_x\in V\}$. Since $\sigma_x$ commutes with
$G$, the subspace $E_0$ is a $G$-submodule. As
$V\le E^{\langle x\rangle}=E\sigma_x$ and $\dim V=1$, we have $\dim E_0=n/2+1$, which is less than $n$ as
$n\ge3$ by Lemma~\ref{local}. Since $G$ is a $2$-group, the nonzero
module $E/E_0$ has a trivial quotient of dimension one. Hence there is
a nonzero $G$-equivariant linear map $s:E\to V$ vanishing on $E_0$. It
vanishes on $V$ as $V\le E_0$, and on $W$ as $G$ acts trivially on
$V$.

By Lemma~\ref{lift-transvection}, there is an automorphism
$\alpha:g\mapsto g\delta(g)$ of $P$ inducing the identity on $P/E$
with $\alpha|_E=1+s$. It fixes $W$ pointwise, as $s(W)=0$. By
Lemma~\ref{transfer}\textup{(iv)}, we may assume that
$\delta(g)\sigma_x=\delta(g)(x-1)$ lies in $V$ for all $g\in P$, that
is, $\delta(g)\in E_0\le\ker s$. By Lemma~\ref{iterate}, $\alpha^2(g)=g\,s(\delta(g))=g$. Thus $\alpha$
is an involution. It fixes $[P,\alpha]$ pointwise, since $[P,\alpha]$
is generated by the elements $\delta(g)\in\ker s$. Hence $\alpha$ is
quadratic, and it is noninner by Lemma~\ref{outer}. This completes the
proof.
\end{proof}

\subsection{Detecting the extension class}
By Proposition~\ref{four}, we may assume that $G$ has no normal
four-subgroup. Then every normal abelian subgroup of $G$ is cyclic,
since a noncyclic normal abelian subgroup $A$ of $G$ would give a
normal subgroup of order four of $G$ inside the elementary abelian
group $\Omega_1(A)$. By \cite[Theorem~5.4.10(i)]{Gor}, the noncyclic
group $G$ is dihedral of order at least sixteen, semidihedral of order
at least sixteen, or generalized quaternion. In each case $Z(G)$ has order two, and so $U=Z(G)$. If $R$
is a cyclic subgroup of index two in $G$, then the unique involution of
$R$ is central in $G$, since $R$ is normal. Hence $U\le R$.

\begin{lemma}\label{support}
Assume that $G$ is dihedral of order at least sixteen, or that $G\cong Q_8$, or that
$\dim E^R=2$ for some cyclic subgroup $R$ of index two in $G$. Then
$P$ has a quadratic noninner involution with displacement in $V$.
\end{lemma}
\begin{proof}
By Lemma~\ref{supplement}, it suffices to find a proper subgroup $H$
of $P$ with $P=HE$.

First assume that $G$ is dihedral, say of order $2m\ge16$, and let $R$
be its cyclic subgroup of index two. Choose involutions
$x_1,x_2\in G\setminus R$ whose product generates $R$. They are not
central in $G$, so they do not lie in $S$, whose only involution $u$
is central. By Lemma~\ref{local}, they lift to involutions
$t_1,t_2\in P$. Put $p_0=t_1t_2$. Then $p_0^{m/2}\in N\setminus E$,
and so $p_0^m=v$ by Lemma~\ref{local}. Hence $p_0$ has order $2m$, and
$H=\langle t_1,t_2\rangle$ is a dihedral group of order $4m$. Since
$H$ maps onto $G$, we have $P=HE$ and $|H\cap E|=2$. As $n\ge3$, the
subgroup $H$ is proper, as desired.

Next, let $R$ be any subgroup of index two in $G$. Then $G/R$ acts on
$E^R$ with fixed space $E^G=V$, and so $\dim E^R\le2$. Assume that $R$
is cyclic and $\dim E^R=2$. Since $G/R$ fixes only the vectors of $V$,
it interchanges the two vectors of $E^R\setminus V$. These two vectors
form a basis of $E^R$, and so $E^R\cong\mathbb F_2[G/R]$. Let
$\kappa\in H^2(G,E)$ be the class of the extension $E\to P\to G$.
Since $U$ acts trivially on $E$, we have $H^2(U,E)\cong E$, and the
restriction of $\kappa$ to $U$ corresponds to the square of a lift of
$u$, which is $v$ by Lemma~\ref{local}. Let $c\in H^2(G,E^R)$
correspond under Shapiro's lemma to the class of the cyclic double
cover of $R$. Restriction to $R$ followed by evaluation at the
identity coset is the Shapiro map
\cite[III.5.9 and III.8, Exercise~2]{Brown}. Since $U\le R$ and the cover does
not split over $U$, evaluation sends the restriction of $c$ to $U$ to
a nonzero element. On the other hand, $U$ is central in $G$ and acts
trivially on $E^R$, so $G$ acts on $H^2(U,E^R)\cong E^R$ through its
action on $E^R$. Restrictions from $G$ are invariant under this action:
in an extension representing $c$, a lift $\widetilde u$ of $u$
satisfies $\widetilde u^{\widetilde g}=\widetilde u e'$ with $e'\in E^R$
for every lift $\widetilde g$ of $g\in G$.
Since $u$ acts trivially on $E^R$, we have
$(\widetilde u^2)^{\widetilde g}=(\widetilde u e')^2=\widetilde u^2$.
Thus the restriction of $c$ to $U$ lies in $(E^R)^G=V$. Hence it
corresponds to $v$. Thus the image of $c$ in $H^2(G,E)$ and $\kappa$
have the same restriction to $U$, and so they are equal by
Lemma~\ref{local}. It follows that the image of $\kappa$ in
$H^2(G,E/E^R)$ is zero, that is, the extension $E/E^R\to P/E^R\to G$
splits. A complement has the form $H/E^R$, where $P=HE$ and
$H\cap E=E^R$. Since $E^R<E$ as $n\ge3$, the subgroup $H$ is proper.

Finally, assume that $G\cong Q_8$. Since $E/V\ne0$, there is a line in
$(E/V)^G$. Let $E_1$ be its inverse image in $E$. Then $\dim E_1=2$,
and $G$ acts nontrivially on $E_1$ as $E^G=V$. The kernel of this
action is a subgroup $R$ of index two in $Q_8$, hence cyclic, and
$E_1\le E^R$. As $\dim E^R\le2$, we get that $\dim E^R=2$, and the
previous case applies. This completes the proof.
\end{proof}

\subsection{The action on a cyclic subgroup of order four}
We are left with the case where $G$ is generalized quaternion or
semidihedral of order at least sixteen and $E^R=V$ for the cyclic
subgroup $R$ of index two in $G$. We use the presentation
$$
 G=\langle r,t\mid r^m=1,\quad r^t=r^k,\quad t^2=r^l\rangle,
 \qquad R=\langle r\rangle,\quad m\ge8,
$$
where $(k,l)=(-1,m/2)$ in the generalized quaternion case and
$(k,l)=(m/2-1,0)$ in the semidihedral case. In both cases
$k\equiv-1\pmod4$, as $m\ge8$, and $U=\langle r^{m/2}\rangle$. Since
$Y/E$ is a nontrivial central subgroup of $G$, we get that $Y/E=U$.
As $Y/E\cong Y^2=Z$, this gives $Z=V$.

Put $X=r-1$, acting on $E$ from the right. Since $E^R=V$ is a line,
$r$ acts on $E$ as a single unipotent Jordan block of size $n$. As
$r^{m/2}=u$ acts trivially on $E$, we have $X^{m/2}=r^{m/2}-1=0$ on
$E$, and so $n\le m/2$. Recall that $n\ge3$ by Lemma~\ref{local}.
Since $r^k$ generates $R$, each of $r-1$ and $r^k-1$ is a multiple of
the other in $\mathbb F_2R$. As $t$ acts on $R$ by $r\mapsto r^k$, each
subspace $EX^i$ is $G$-invariant. Each layer $EX^i/EX^{i+1}$ has
dimension one, so $G$ acts trivially on it. In particular, $EX$ is a
$G$-invariant hyperplane with trivial quotient, and so $EX=W$.

Choose lifts $p_0,h\in P$ of $r,t$, and put $b=p_0^{m/2}$. Then
$bE=u$, and $b^2=v$ by Lemma~\ref{local}. We claim that
\begin{equation}\label{conj-b}
 b^{p_0}=b,\qquad b^h=b^{-1}.
\end{equation}
Since $b$ is a power of $p_0$ and $b\in D$, it centralizes
$\langle E,p_0\rangle$, which has index two in $P$. As $b\notin E$ and
$Z=V\le E$, the element $b$ is not central, and so its conjugates are
precisely $b$ and $b^h$. Put $e=b^{-1}b^h$. Then $1\ne e\in E$, since
$bE=u$ is central in $G$. As $b,b^h\in N\le D$ commute,
conjugation permutes $b$ and $b^h$, and
hence sends $e$ to $e$ or $e^{-1}=e$. Thus $e\in Z=\langle b^2\rangle$,
and so $e=b^2$. It follows that $b^h=b^{-1}$, as desired.

Let $s:E\to V$ be the linear map with kernel $W$. It is
$G$-equivariant and vanishes on $V$. By Lemma~\ref{lift-transvection},
there is an automorphism $\alpha:g\mapsto g\delta(g)$ of $P$ inducing
the identity on $P/E$ with $\alpha|_E=1+s$. By Lemma~\ref{iterate},
$\alpha^2(g)=g\,s(\delta(g))$. Since $\delta$ is a derivation and $G$
acts trivially on $V$, the map $g\mapsto s(\delta(g))$ is a
homomorphism $P\to V$. It vanishes on $E$ as $s^2=0$. Let
$\theta:G\to\mathbb F_2$ be the character with
$s(\delta(g))=v^{\theta(gE)}$.

By the derivation identity,
$\delta(b)=\delta(p_0)(1+r+\cdots+r^{m/2-1})$. As $m/2$ is a power of
two, we have $1+r+\cdots+r^{m/2-1}=X^{m/2-1}$. Thus $\delta(b)=0$ if
$n<m/2$. We treat the case $n=m/2$ in Proposition~\ref{full}.

\begin{lemma}\label{correction}
For every character $\lambda:G\to\mathbb F_2$, there is a derivation
$f:G\to\langle b\rangle$ with $f(g)E=u^{\lambda(g)}$ for all $g\in G$
and $f(u)=1$.
\end{lemma}
\begin{proof}
By \eqref{conj-b}, the subgroup $\langle b\rangle$ is normal in $P$,
and $G$ acts on it with $r$ acting trivially and $t$ inverting it.
Let $\langle b\rangle\rtimes\langle j\rangle$ be the semidirect product
in which the involution $j$ inverts $b$. We claim that
$r\mapsto b^{\lambda(r)}$ and $t\mapsto jb^{\lambda(t)}$ define a
homomorphism from $G$. To see this, we check the defining relations
of $G$. First, $(b^{\lambda(r)})^m=1$ as $4\mid m$.
Since $k\equiv-1\pmod4$, we have
$$
 (jb^{\lambda(t)})^{-1}b^{\lambda(r)}jb^{\lambda(t)}=b^{-\lambda(r)}
 =b^{k\lambda(r)}.
$$
Finally, $(jb^{\lambda(t)})^2=1=b^{l\lambda(r)}$ since $4\mid l$. Let
$\omega:G\to\mathbb F_2$ be the character with kernel $R$. We write the
image of $g$ as $j^{\omega(g)}f(g)$ with $f(g)\in\langle b\rangle$. The
multiplication rule in the semidirect product shows that
$f(gh)=f(g)^hf(h)$, and $f$ modulo $\langle b^2\rangle$ is given by
$\lambda$. Thus $f(g)E=u^{\lambda(g)}$. Finally,
$f(u)=f(r^{m/2})=b^{\lambda(r)m/2}=1$, since $4\mid m/2$. This
completes the proof.
\end{proof}

\begin{proposition}\label{corrected}
If $n<m/2$, then $P$ has a cubic noninner involution $\beta$ with
$[P,\beta]\le N$.
\end{proposition}
\begin{proof}
Let $f$ be as in Lemma~\ref{correction} for $\lambda=\theta$, and put
$$
 \beta(g)=g\,\delta(g)f(gE).
$$
Both $\delta$ and $g\mapsto f(gE)$ are derivations with values in the
abelian normal subgroup $N=E\langle b\rangle$, and hence so is their
product. Thus $\beta$ is an endomorphism of $P$. It acts as $1+s$ on
$E$, and it fixes $b$ since $\delta(b)=0$ and $f(u)=1$. Let $g\in P$,
and write $f(gE)=b^j$. Then $j\equiv\theta(gE)\pmod2$, and so
$$
 \beta^2(g)=g\,s(\delta(g))\,b^{2j}=g\,v^{\theta(gE)+j}=g.
$$
Thus $\beta^2=1$, so $\beta$ is bijective, and $\beta$ is an involution
as $\beta\ne1$. Its successive displacements lie in
$N$, $V$ and $1$, so $\beta$ is cubic and $[P,\beta]\le N$. It is
noninner by Lemma~\ref{outer}. This completes the
proof.
\end{proof}

\begin{proposition}\label{full}
If $n=m/2$, then the central extension $E/W\to P/W\to G$ splits, and
$P$ has a quadratic noninner involution with displacement in $V$.
\end{proposition}
\begin{proof}
Put $\sigma=1+r+\cdots+r^{m/2-1}$. Then $\sigma=X^{m/2-1}=X^{n-1}$ on
$E$, and its kernel is $EX$. Write $p_0^h=p_0^ke_1$ with $e_1\in E$.
Since $k$ is odd and $r^{m/2}$ acts trivially on $E$, we also have
$\sigma=1+r^k+\cdots+r^{k(m/2-1)}$ on $E$. Raising $p_0^h=p_0^ke_1$ to
the power $m/2$ and using \eqref{conj-b}, we obtain that
$$
 b^{-1}=b^h=b^k\cdot e_1\sigma=b^{-1}\cdot e_1\sigma,
$$
where $e_1\sigma$ is written multiplicatively. Hence $e_1\sigma=0$,
and so $e_1\in EX=W$.

Since $r-1$ and $r^k-1$ are multiples of each other in
$\mathbb F_2R$, we have $E(r^k-1)=EX$. Hence we may replace $h$ by
$he$ for a suitable $e\in E$ so that $p_0^h=p_0^k$. Write
$h^2=p_0^le_0$ with $e_0\in E$. Conjugating $p_0$ by $h$ twice gives
$$
 p_0[p_0,e_0]=p_0^{h^2}=p_0^{k^2}.
$$
As $k^2\equiv1\pmod m$ and $p_0^m=b^2=v$, we get that
$[p_0,e_0]\in V$. In additive notation, $[p_0,e_0]=e_0X$, and so
$e_0X\in V=\ker X$. Since $X$ is a single Jordan block of size
$n\ge3$, we have $\ker X^2=EX^{n-2}\le EX$. Thus $e_0\in EX=W$.

Consequently, in $P/W$, the images of $p_0$ and $h$ satisfy the
defining relations of $G$, including $p_0^m=1$. Hence
$H=\langle p_0,h\rangle W$ satisfies $P=HE$ and $H\cap E=W$, and so
the extension splits. Since $W<E$, the subgroup $H$ is proper, and
Lemma~\ref{supplement} completes the proof.
\end{proof}

\begin{proof}[Proof of Theorem~\ref{two}]
Let $E$ be maximal among the normal elementary abelian subgroups of
$P$. If $E$ is unsaturated, a noninner element of $D_E(P)$ is a
quadratic noninner involution. Thus we may assume that $E$ is
saturated. If $b_G(E)>0$, Proposition~\ref{invertcase} applies. If
$D=E$, Proposition~\ref{split-two} applies. In the remaining case
$b_G(E)\le0$ and $D>E$. If $G$ has a normal four-subgroup,
Proposition~\ref{four} applies. Otherwise $G$ is dihedral,
semidihedral or generalized quaternion. Then Lemma~\ref{support} and
Propositions~\ref{corrected} and~\ref{full} cover all cases, since
$n\le m/2$. The second assertion was proved in
Section~\ref{sec:two}. This completes the proof.
\end{proof}

\section{The subgroup \texorpdfstring{$K_P(E)$ for $p=2$}{K\_P(E) for p=2}}\label{sec:Kp}

For odd $p$, the automorphism in Theorem~\ref{flag} belongs to
$D_F(P)$, and $K_P(F)=C_P(E)$ by Section~\ref{sec:flag}. In this
section, we show that for $p=2$ the automorphism in
Theorem~\ref{two} can be chosen with displacement in $K_P(E)$, unless
$P\cong D_8$.

Let $Y(B)$, $C(B)$ and $K_P(E)$ be as in Section~\ref{sec:flag}. Thus
$C(B)$ consists of the elements $x\in C_P(B)$ with $[P,x]\le B$ and
$x^2\in B$, and $Y(E)=Y$. Put $K=K_P(E)$, and let $M$ be the subgroup
at which the series defining $K$ stabilizes. Then $Z(K)=M$ and
$K/M\le\Omega_1(Z(P/M))$. Since $E\le Y$ and $Y$ centralizes $E$, we
also have
$$
 E\le Z(Y)=M_0\le M\le K.
$$

For a characteristic abelian subgroup $U$, the following properties are part
of \cite[Theorem~D and Corollary~E]{Kizmaz}. Here $Z(Y)$ need not be
characteristic, so we include the short argument.

We claim that $C_P(K)=M$. Clearly $C_P(K)\cap K=Z(K)=M$. Assume that
$C_P(K)>M$. Then the nontrivial normal subgroup $C_P(K)/M$ of $P/M$
contains a central involution of $P/M$. Any lift of it centralizes
$M$, and so it belongs to $C(M)=K$. Thus it lies in $C_P(K)\cap K=M$,
which is a contradiction. Together with $K/M\le\Omega_1(Z(P/M))$, this
shows that $K$ satisfies conditions (i)--(iii) of Thompson's critical
subgroup theorem \cite[Theorem~5.3.11]{Gor}, although it need not be
characteristic in $P$.

\begin{proposition}\label{located}
Let $P$ be a finite nonabelian $2$-group, and let $E$ be maximal among
its normal elementary abelian subgroups. If $P\not\cong D_8$, then the
automorphism $\alpha$ in both assertions of Theorem~\ref{two} can be
chosen with $[P,\alpha]\le K_P(E)$. If $P\cong D_8$, then no noninner
automorphism $\alpha$ of $P$ satisfies $[P,\alpha]\le K_P(E)$.
\end{proposition}
\begin{proof}
We follow the proof of Theorem~\ref{two}. If $E$ is unsaturated, a
noninner element $\alpha$ of $D_E(P)$ satisfies $[P,\alpha]\le E\le K$.
Thus we may assume that $E$ is saturated.

First assume that $b_G(E)>0$. By Proposition~\ref{large}, we have
$Y=EJ=D$ and $Z(Y)=EZ=E$. Every element of $Y$ centralizes $E$, has its
square in $Z\le E$, and commutes with every element of $P$ modulo $E$.
Thus $C(E)=Y$. Since $M_0=Z(Y)=E$, the series stabilizes at $M=E$, and
so $K=Y$. If $Y=P$, there is nothing to prove. Otherwise, let
$c_t\gamma$ be the involution constructed in the proof of
Proposition~\ref{invertcase}. The automorphism $\gamma$ has
displacement in $B=E\langle y\rangle\le Y$, and $c_t$ has displacement
in $E$ as $t\in Y$. Hence $[P,c_t\gamma]\le Y=K$.

Next assume that $D=E$. Since $P\not\cong D_8$,
Proposition~\ref{split-two} gives a quadratic noninner involution with
displacement in $V\le E\le K$.

Now assume that $b_G(E)\le0$ and $D>E$. Then $K=D$. To see this, note
that $D$ is abelian by Proposition~\ref{structure}. Hence
$K=C(M)\le C_P(M)\le C_P(E)=D$, and so $D\le C_P(K)=M\le K$. The
involutions constructed in Proposition~\ref{four}, Lemma~\ref{support}
and Proposition~\ref{full} have displacement in $E$, and the involution
of Proposition~\ref{corrected} has displacement in $N\le D$. Thus all of
them have displacement in $K$.

In the second assertion of Theorem~\ref{two}, the involution comes from
the unsaturated case or from Proposition~\ref{split-two}, and so it has
displacement in $K$ as well.

Finally, let $P\cong D_8$, written as in Example~\ref{d8}. By symmetry,
we may assume that $E=\langle r^2,t\rangle$. Then $C_P(E)=E$, and so
$Y(E)=E$, $M_0=E$ and $C(E)=E$. Hence $K_P(E)=E$. If $[P,\alpha]\le E$
for an automorphism $\alpha$, then $\alpha$ maps $E$ to itself.
However, by Example~\ref{d8}, every noninner automorphism of $D_8$
sends $t$ to $r^jt$ with $j$ odd, and so it interchanges the two normal
four-subgroups $\langle r^2,t\rangle$ and $\langle r^2,rt\rangle$. This
completes the proof.
\end{proof}

\section*{Acknowledgments}
The author used large language models to explore proof strategies,
examine individual arguments, revise the exposition, search the
literature, and prepare GAP computations.
The author takes responsibility for the arguments and references
in this paper.

\end{document}